\documentclass{article}
\usepackage{amssymb,amsmath,amsfonts,enumerate,url,graphicx}
\usepackage{enumerate}
\usepackage{subfigure}

\newtheorem{theorem}{Theorem}[section]

\newtheorem{definition}[theorem]{Definition}

\begin{document}

\title{Exploration of offsets of Cayley ovals and their singularities}

\author{Thierry Dana-Picard and Daniel Tsirkin}

\maketitle

\begin{center}
Department  of  Mathematics,  Jerusalem  College  of
Technology\\  
Jerusalem 9116011, Israel\\e-mail: ndp@jct.ac.il, academictsirkin@gmail.com
\end{center}



\begin{abstract}
We explore offsets of Cayley ovals, by networking with different kinds of software. Using their specific abilities, algebraic, geometric, dynamic, we conjecture interesting properties of the offsets. For a given progenitor (the given plane curve whose offsets are studied), changes in the offset distance induce great changes in the shape and the topology of the offset.  Such a study has been performed in the past for classical curves, and recently for non classical ones.Here we relate to Cayley ovals; despite them being non singular, their offsets have intriguing properties, cusps, and self-intersections. We begin with a short study of envelopes of families of circles with constant radius centered on the oval (these constructs are often studied together with offsets, but they are different objects). Then we study the offsets, which are defined as geometric loci. Both approaches are supported by the automated methods provided by the software.
\end{abstract}

\maketitle

\section{Introduction}

\subsection{Geometric loci in a technology rich environment.}
The study of (real)plane algebraic curves is an active domain of research. Among the topics under study are the determination of geometric loci \cite{botana-abanades,bp}, envelopes \cite{botana recio,dpz - revival}, isoptic curves \cite{DZM bisoptics of hyperbolas,dpk dynamics,DPM-inner} etc. Recently AI systems entered the environment and their affordances are under scrutiny \cite{emprin-richard,BDR,botana recio AMAI}; we think that the generative AI systems widely available are still not ripe for mathematical work, therefore, we do not use them here.

Algorithms for automated methods in Geometry have been developed and implemented in various kinds of software. We use GeoGebra  for its strong features for dynamic exploration (e.g.,  dragging points, slider bars, etc.); see \cite{ART-ecosystem,botana-valcarce 2}. A given curve can be build via different methods, making the study richer. A Computer Algebra System called Giac is also implemented in GeoGebra \cite{kovacs parisse}. but sometimes stronger abilities are needed for polynomials of higher degree. For most algebraic manipulations, we use the Maple 2024 Computer Algebra System (CAS). Networking between the two kinds of software is central for a thorough exploration of the curves under study \cite{DPK dialog}. Despite wishes expressed for a long time of more automatic dialog between different kinds of software \cite{roanes 2003,roanes 2020}, the data transfer from one to the other has still to be done by hand.   

\subsection{Envelopes and offsets}

Bruce and Giblin \cite{bruce and giblin} give 4 non-equivalent definitions of an envelope of a 1-parameter family of plane curves.  In this paper we work according to one of them, which is the only definition of an envelope mentioned by Berger \cite{berger}. Kock \cite{kock} calls it \emph{analytical definition} of an envelope.

\begin{definition}
Let $\mathcal{C}_t$ be a 1-parameter family of curves given by an equation $F(x, y, t) = 0$, where $t$ is a real parameter. An envelope of the family $\mathcal{C}_t$ is the set of all points $(x, y)$ in the plane verifying the system of equations
\begin{equation}
\label{eq envelope}
\begin{cases}
F(x, y, t) = 0 \\
\frac{\partial}{\partial t}F(x, y, t) = 0
\end{cases}.
\end{equation}
\end{definition}                                                              

\begin{figure}[htb]
\begin{center}
\includegraphics[width=4.5cm]{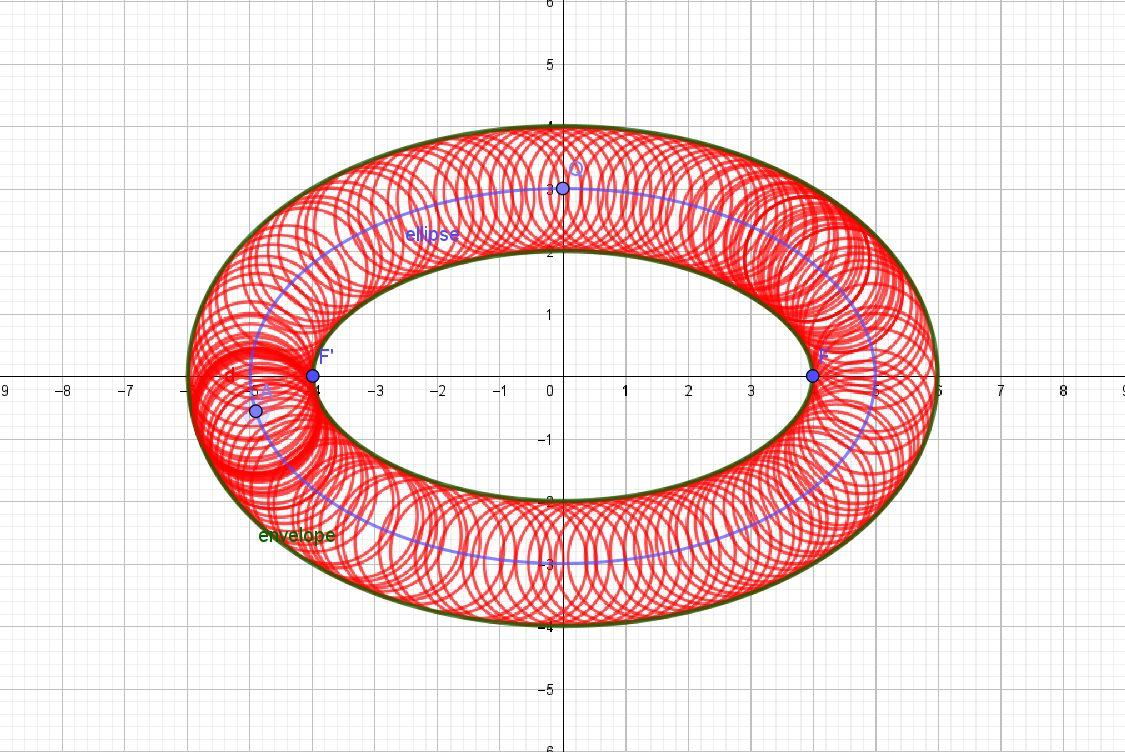}
\qquad
\includegraphics[width=4.5cm]{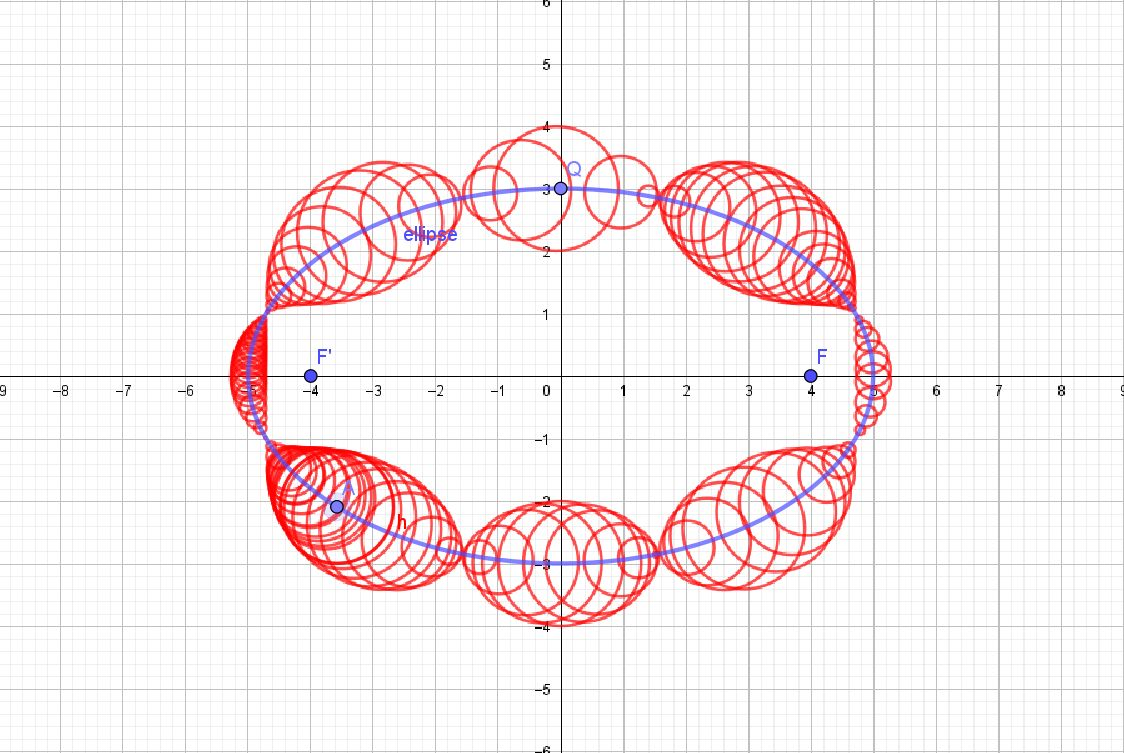}
\caption{Envelopes of two families of circles centred in the same ellipse}
\label{fig 1 envelopes}
\end{center}
\end{figure}

Figure \ref{fig 1 envelopes} illustrates the construction of an envelope of (a) the family of unit circles centered on the ellipse whose equation is  $\frac{x^2}{25}+\frac{y^2}{9}=1$, (b) the family of circles centered on the same ellipse and whose radius is $\vert \cos x_A \vert$, when $A$ is the center of the circle. The pictures have been obtained from animations with GeoGebra\footnote{They are  available at (a) \url{https://www.geogebra.org/m/e4ahaq8r}  and  (b) \url{https://www.geogebra.org/m/qfnh7czy}.}    

Using GD \textbf{Envelope} command, we obtained also the following equation for the envelope: 
\footnotesize
\begin{align*}
& 81x^8 + 612x^6 y^2 - 6516x^6 + 1606x^4 y^4 - 32444x^4 y^2 + 182848x^4 \\
& + 1700x^2 y^6 - 35612x^2 y^4 + 468352x^2 y^2 - 2179072x^2 + 625y^8 + 6700y^6 \\
& - 196544y^4 - 1720320y^2 + 9437184=0.
\end{align*}
\normalsize

Using a CAS, it can be proven that the polynomial in the left hand side is irreducible, whence the impossibility to distinguish the two loops of the envelope by algebraic computations.

The construction of the envelope of circles with variable radius involves steps which are not supported by the algorithm of the \textbf{Envelope} command, therefore an automated equation for the envelope cannot be obtained.

\subsubsection{Offsets}
\label{subsub offsets}
Another object, sometimes identified with an envelope, is the so-called offset of a curve $\mathcal{C}$; Leibniz \cite{leibnitz} called this a parallel curve to $\mathcal{C}$.
\begin{definition}
\label{def offset}
Let $\mathcal{C}$ be a given plane curve. At each regular point $A$ of $\mathcal{C}$, denote by $\overset{\longrightarrow}{n_A}$  a unit normal vector to $\mathcal{C}$. Let $d$ be a positive real number. Now denote by $B$ a point in the plane such that $\overset{\longrightarrow}{AB}=d  \overset{\longrightarrow}{n_A}$ . The geometric locus of the point $B$ when $A$ runs over the curve $\mathcal{C}$ is called an offset of $\mathcal{C}$ at distance $d$. The curve $\mathcal{C}$ is called the \emph{progenitor} of the offset.
\end{definition}

Note that the vector  $\overset{\longrightarrow}{n_A}$ is not uniquely defined; therefore, there may exist two offsets at the same distance, we will call them internal and external offsets (in general this make sense in an obvious way). Actually, we will mention "the" offset at distance $d$ as the union of these components. See Section \ref{section offsets of Cayley ovals}.

\begin{figure}[htb]
\begin{center}
 \includegraphics[width=5.5cm]{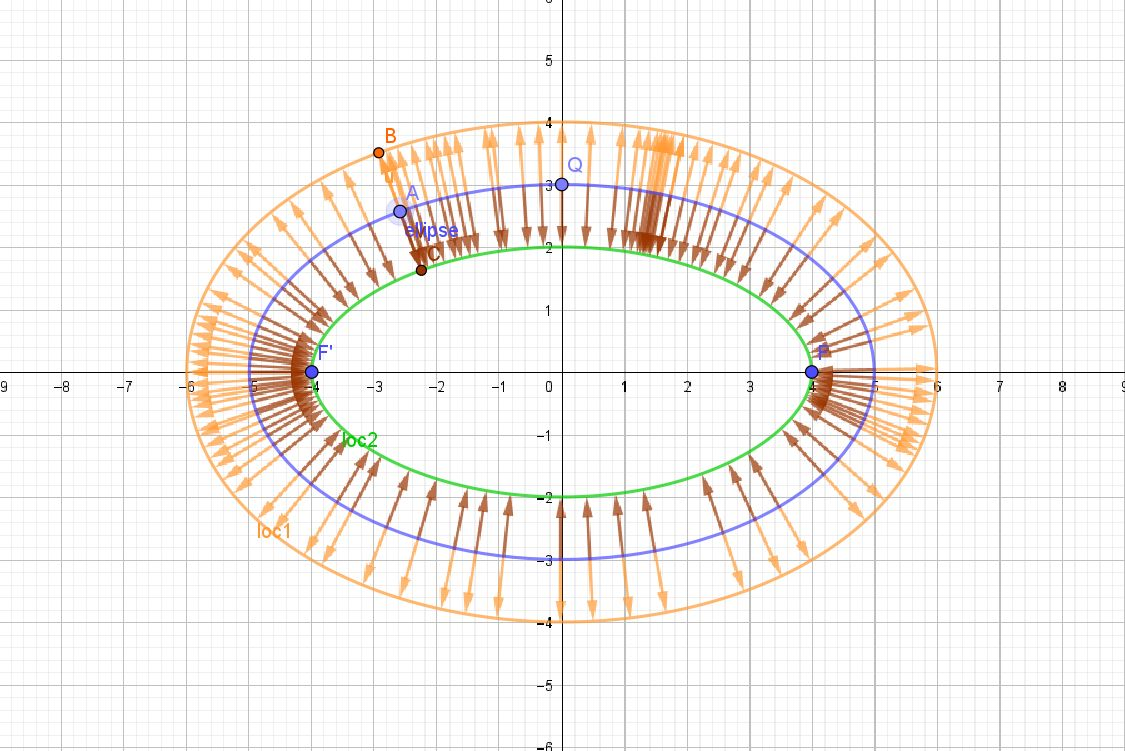}
 \qquad
 \includegraphics[width=5.5cm]{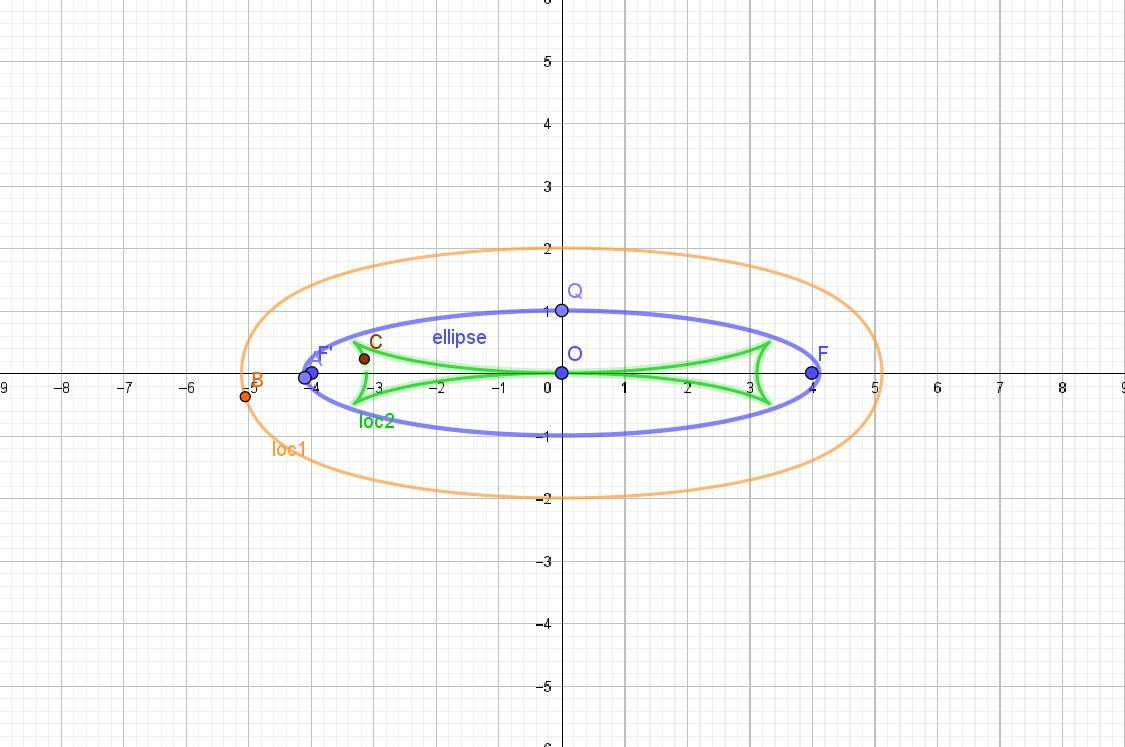}
\caption{Offsets at distance 1 of two different ellipses}
\label{fig 2 offsets}
\end{center}
\end{figure}

Figure \ref{fig 2 offsets}(a) shows the offset at distance 1 of an ellipse; this offset is also the envelope of the family of unit circles centered on the ellipse. The picture is a screenshot of a GeoGebra session\footnote{\url{https://www.geogebra.org/m/sskh84rd }}; the construction is as follows: 
\begin{itemize}
\item For a point $A$ on the ellipse, GeoGebra's pre-defined commands (available as buttons and as written commands) provide a unit circle centered at $A$, the tangent to the ellipse at $A$, the normal to the ellipse at $A$ and the intersection of this normal with the pre-obtained circle. These are the points $B$ and  $C$, whose geometric locus is the desired offset.  
\item The envelope is previewed by moving the point $A$ on the ellipse with \emph{Trace On} for the circle. Some steps are not supported by the \textbf{Envelope} command. 
item The offset at distance 1 is previewed viewed by moving the point $A$ on the ellipse with \emph{Trace On} for points $B$ and  $C$. Moreover it can be obtained also with the automated command \textbf{Locus(Path, Point)}. 
Note that the obtained curves seem to coincide, but we do not have here an exact proof of this coincidence, as the commands are numerical.  
\end{itemize}
Figure \ref{fig 2 offsets} shows the construction of offsets at distance 1 of two different ellipses. The first one has equation  $\frac{x^2}{16}+\frac{y^2}{9}=1$, and the second one has equation $\frac{x^2}{17}+y^2=1$
 Note that this time, the external component seems non-singular, but the internal component has cusps and, maybe, a double point. This is a good example of something well-known: the topology of an offset may be must more complicated than the topology of the progenitor; see \cite{alcazar sendra 2007}.

\subsubsection{Identification of components}
\begin{figure}[htb]
\begin{center}
 \includegraphics[width=4.5cm]{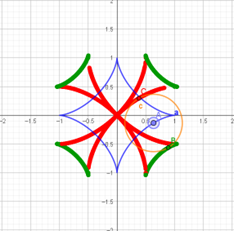}
 \qquad
 \includegraphics[width=4.5cm]{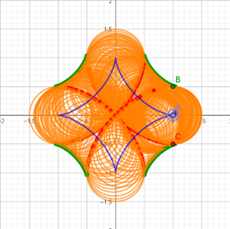}
\caption{Comparison of envelope and offset of a singular curve}
\label{fig envelope-offset}
\end{center}
\end{figure}

 When the progenitor has singularities, the envelope of circles centred on it with radius $d$ and the offset at distance $d$ may be very different. Such occurrences have been shown in \cite{offsets astroid}, from which Figure \ref{fig envelope-offset} is taken\footnote{The GeoGebra applet is available at \url{https://www.geogebra.org/m/ucuvdwqe}.}. 

GeoGebra-Discovery's  \textbf{LocusEquation} provides a polynomial equation of degree 23 for the geometric locus. Let us denote it by $P(x,y)$.  The CAS implemented in  GeoGebra provides a factorization of  $P(x,y)$, but hardware constraints led us to perform the same task with Maple. The following result is obtained:
 \begin{enumerate}[(i)]
\item Three linear factors: $F_1(x,y)=x$, $F_2(x,y)=x-1$, $F_3(x,y)=x+1$.
\item One quartic factor: $F_4(x,y)=x^{4}+2 x^{2} y^{2}+y^{4}-34 x^{2}+30 y^{2}+289$
\item Two octic factors: 
\begin{align*}
F_5(x,y) & =x^8+36 x^6 y^2+358 x^4 y^4+612 x^2 y^6+289 y^8-68 x^6-1772 x^4 y^2\\
\quad & \qquad -9996 x^2 y^4+8092 y^6+1700 x^4+27336 x^2 y^2+47396 y^4\\
\quad & \qquad -18496 x^2-129472 y^2+73984\\
F_6(x,y)& =x^{8}+36 x^{6} y^{2}+358 x^{4} y^{4}+612 x^{2} y^{6}+289 y^{8}-36 x^{6}-2252 x^{4} y^{2}\\
\quad & \qquad -11052 x^{2} y^{4}+7548 y^{6}+256 x^{4}+46080 x^{2} y^{2}\\
\quad & \qquad +30720 y^{4}-262144 y^{2}.
\end{align*}
 \end{enumerate}
 The  linear factors are irrelevant and $F_4(x,y)=0$ corresponds to an empty set. The equation $F_6(x,y)=0$ determines the envelope observed with GD, and $F_5(x,y)=0$ is plotted by Maple, but is irrelevant to the present question. It appears for Zariski topological reasons (see \cite{cox,isoptics of quartics}, beyond the scope of the present work.    
At least, what happens here gives a partial explanation why to network between two different kinds of software.

\section{Cayley ovals}
\label{section Cayley ovals}

Ovals are interesting curves, showing numerous geometrical properties and offering a great variety of possible constructions \cite{mazotti}.  
Envelopes of families of circles centered on conics are studied in \cite{dpz - revival}, and centered on Cassini ovals in  \cite{offsets of Cassini ovals}. Here we perform a similar study, making a transition from curves of degree 2, 4 or 6 towards Cayley ovals, which are plane curves of degree 8 (octics). This will add some aspects to the study of octics, from a point of view different from \cite{DPR 2023a, DPR 2023b}.

\begin{definition}[bifocal]
\label{def cayley oval}
Let $F_1$ and $F_2$ be two points in the plane and $b$ a positive real number. The geometric locus of points $M$ in the plane such that the harmonic mean of the distance to $F_1$ and $F_2$ is equal to $b$, i.e., such that $\frac{1}{MF_1}+\frac{1}{MF_2}=\frac2b$, is called a \emph{Cayley Oval}. The points $F_1$ and $F_2$ are its foci.
\end{definition}

\begin{figure}[htb]
\begin{center}
 \includegraphics[width=4cm]{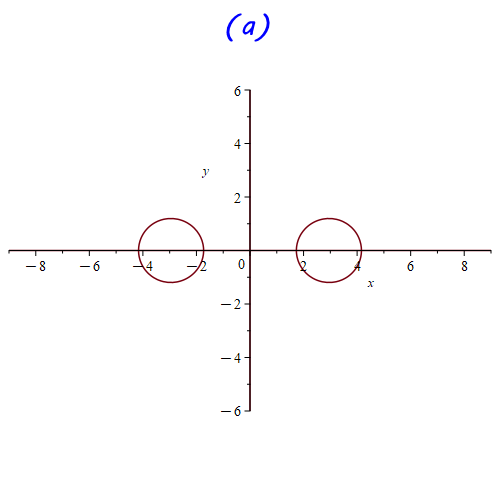}
 \qquad \qquad
 \includegraphics[width=4cm]{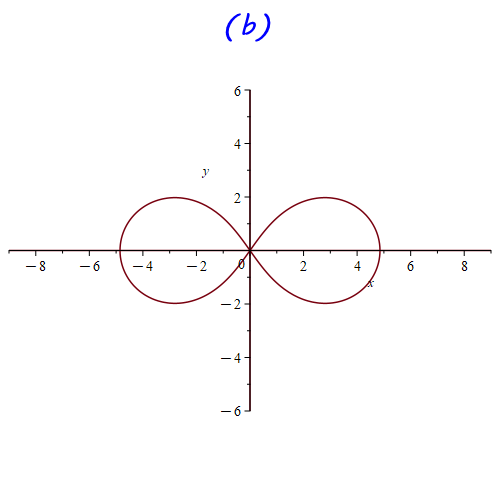}
 \qquad \qquad
 \includegraphics[width=4cm]{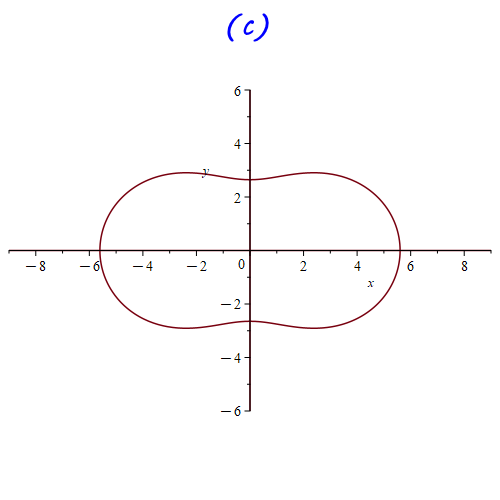}
 \qquad \qquad
 \includegraphics[width=4cm]{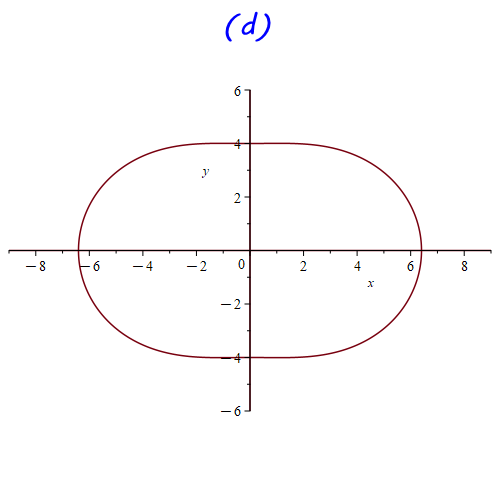}
\caption{Cayley ovals, obtained with the symbolic non polynomial equation.}
\label{fig symbolic non polyn eq}
\end{center}
\end{figure}

Cayley ovals were studied by the 19th English mathematician Arthur Cayley (1821-1895) in 1857. They describe the equipotential lines of the electrostatic potential created by two equal charges placed at the foci (or of the gravitational potential created by two identical masses).  
If $F_1=(a,0)$ and $F_2=(-a,0)$, a Cartesian equation for a Cayley oval is as follows:
\begin{equation}
\label{eq Cayley oval}
\frac{1}{\sqrt{(x+a)^2+y^2}}+\frac{1}{\sqrt{(x-a)^2+y^2}}=\frac{2}{b},
\end{equation}
We use the following Maple code to plot the curves given by Equation (\ref{eq Cayley oval}) and shown in Figure \ref{fig symbolic non polyn eq}.

\small
\begin{verbatim}
F := 1/sqrt((x + a)^2 + y^2) + 1/sqrt((x - a)^2 + y^2) - 2/b;
for b to 6 do
implicitplot({x*y = 0, subs(a = 3, subs(b = 2, F)) = 0}, x = -9 .. 9, y = -6 .. 6);
end do;
\end{verbatim}
\normalsize

As we are interested in the topology of the curves, we can consider this special choice of foci WLOG. With some algebraic manipulations, Equation (\ref{eq Cayley oval}) can be transformed into a polynomial equation of degree 8.
\begin{equation}
\label{param cayley}
  \begin{cases}
  x=\frac{b^2}{a} \cdot \frac{\cos 2t}{\sin^4 2t}\\
  y=\pm \sqrt{\frac{b^2}{4 \cos^4 t}-(x-a)^2}
\end{cases}
\end{equation}

Using the automated commands of GeoGebra-Discovery, in particular the command \textbf{LocusEquation}, we obtain both an equation and a plot. We chose foci at coordinates $(3,0), (-3,0)$.
\begin{figure}[htb]
\begin{center}
 \includegraphics[width=3.5cm]{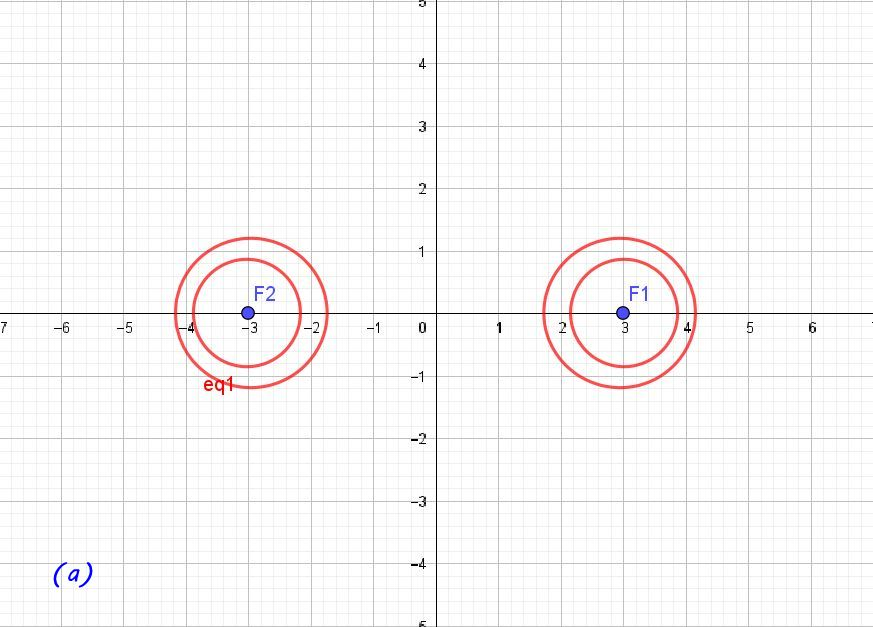}
 \qquad \qquad
 \includegraphics[width=3.5cm]{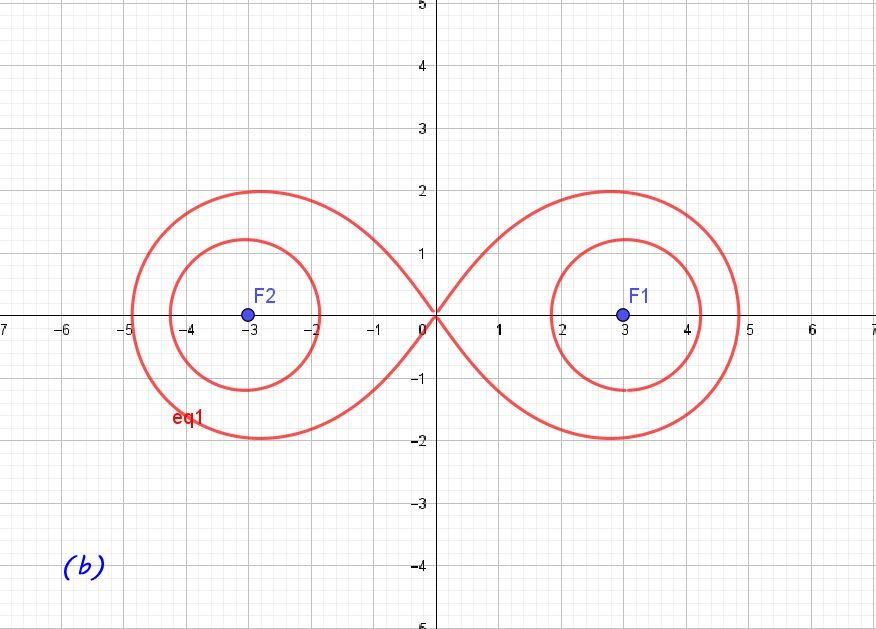}
 \qquad \qquad
 \includegraphics[width=4cm]{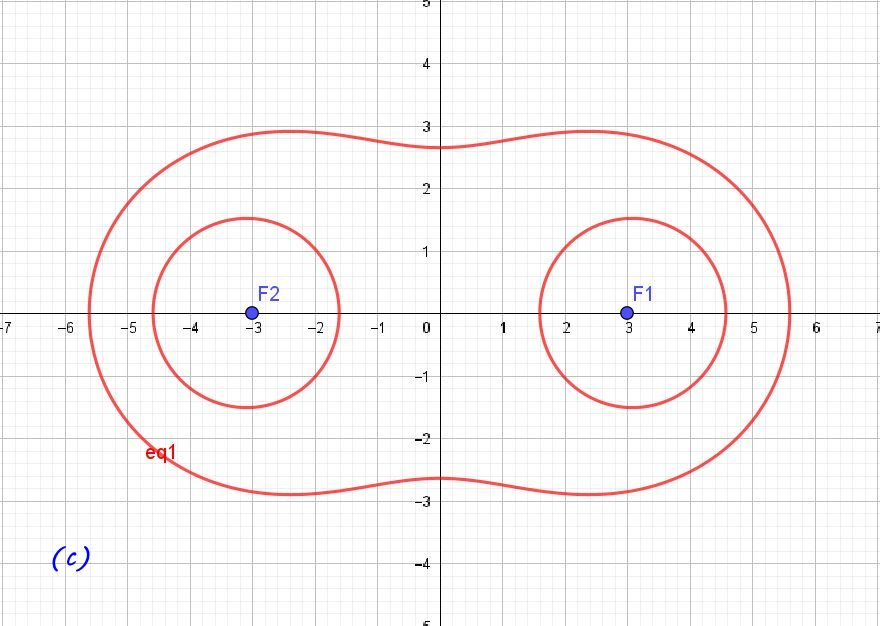}
 \qquad \qquad
 \includegraphics[width=4cm]{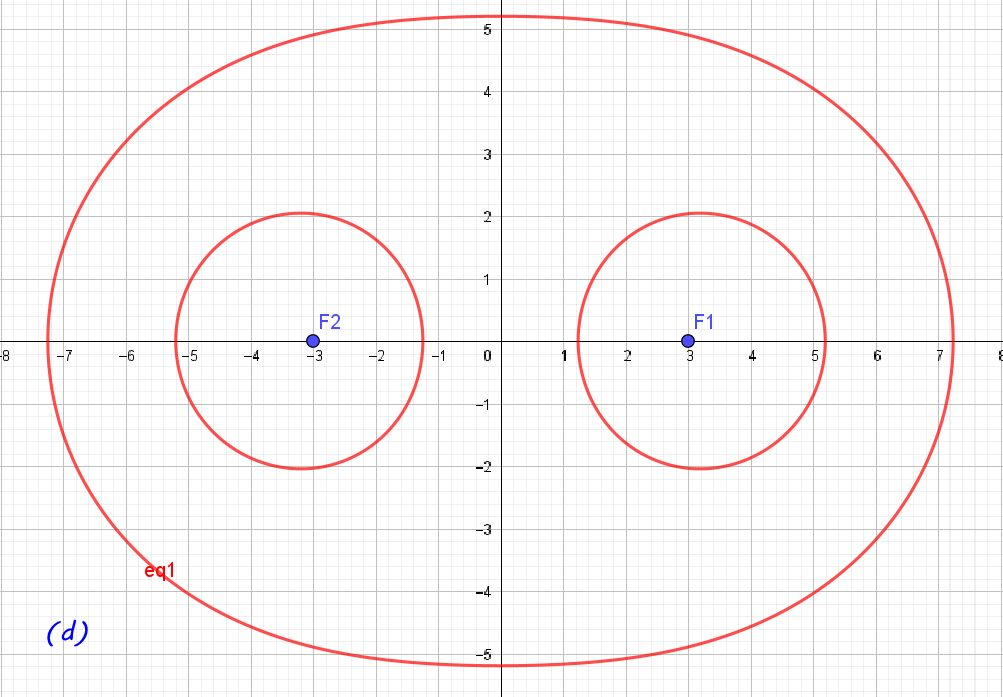}
\caption{Cayley ovals, obtained with GD's LocusEquation}
\label{fig Cayley ovals LocusEquation}
\end{center}
\end{figure}	

The respective equations displayed by the software are: 
\begin{enumerate}[a.]
\item (For $b=2$) $x^8 + 4x^6 y^2 - 40x^6 + 6x^4 y^4 - 48x^4 y^2 + 522x^4 + 4x^2 y^6 + 24x^2 y^4 - 396x^2 y^2 - 2448x^2 + y^8 + 32y^6 + 378y^4 + 1944y^2 + 3645=0$.
\item (For $b=3$; we have a lemniscate) $ x^8 + 4x^6 y^2 - 45x^6 + 6x^4 y^4 - 63x^4 y^2 + 567x^4 + 4x^2 y^6 + 9x^2 y^4 - 486x^2 y^2 - 1458x^2 + y^8 + 27y^6 + 243y^4 + 729y^2 = 0$.
\item (For $b=4$) $x^8 + 4x^6 y^2 - 52x^6 + 6x^4 y^4 - 84x^4 y^2 + 630x^4 + 4x^2 y^6 - 12x^2 y^4 - 612x^2 y^2 + 684x^2 + y^8 + 20y^6 + 54y^4 - 972y^2 - 5103=0$.
\item (For $b=6$) $x^8 + 4x^6 y^2 - 72x^6 + 6x^4 y^4 - 144x^4 y^2 + 810x^4 + 4x^2 y^6 - 72x^2 y^4 - 972x^2 y^2 + 11664x^2 + y^8 - 486y^4 - 5832y^2 - 19683=0$.
\end{enumerate}
Comparing Figures \ref{fig symbolic non polyn eq} and \ref{fig Cayley ovals LocusEquation}, we see that in the last one there are superfluous components (the internal loops in the plots of Figure \ref{fig Cayley ovals LocusEquation}), which do not appear in Figure \ref{fig symbolic non polyn eq} obtained with the original equation (\ref{eq Cayley oval})\footnote{neither with the parametric representation (\ref{param cayley}).}. Why does this occur?

Equation (\ref{eq Cayley oval}) is not a polynomial. In order to derive from it a polynomial equation, algebraic manipulations are needed. They involve twice squaring and isolating the numerator of the expressions. The final output is the following polynomial of degree 8:
\small
\begin{equation}
\label{polyn eq cayley ovals}
\begin{matrix}
P(x,y)= 16 a^{8}-16 a^{6} b^{2}-64 a^{6} x^{2}+64 a^{6} y^{2}+16 a^{4} b^{2} x^{2}-48 a^{4} b^{2} y^{2}+96 a^{4} x^{4}\\
 -64 a^{4} x^{2} y^{2}+96 a^{4} y^{4}+16 a^{2} b^{4} x^{2}+16 a^{2} b^{2} x^{4}-32 a^{2} b^{2} x^{2} y^{2}\\
 -48 a^{2} b^{2} y^{4}-64 a^{2} x^{6}-64 a^{2} x^{4} y^{2}+64 a^{2} x^{2} y^{4}+64 a^{2} y^{6}-16 b^{2} x^{6}\\
 -48 b^{2} x^{4} y^{2}-48 b^{2} x^{2} y^{4}-16 b^{2} y^{6}+16 x^{8}+64 x^{6} y^{2}+96 x^{4} y^{4}+64 x^{2} y^{6}+16 y^{8}=0.
\end{matrix}
\end{equation}
\normalsize

\begin{figure}[htb]
\begin{center}
 \includegraphics[width=7cm]{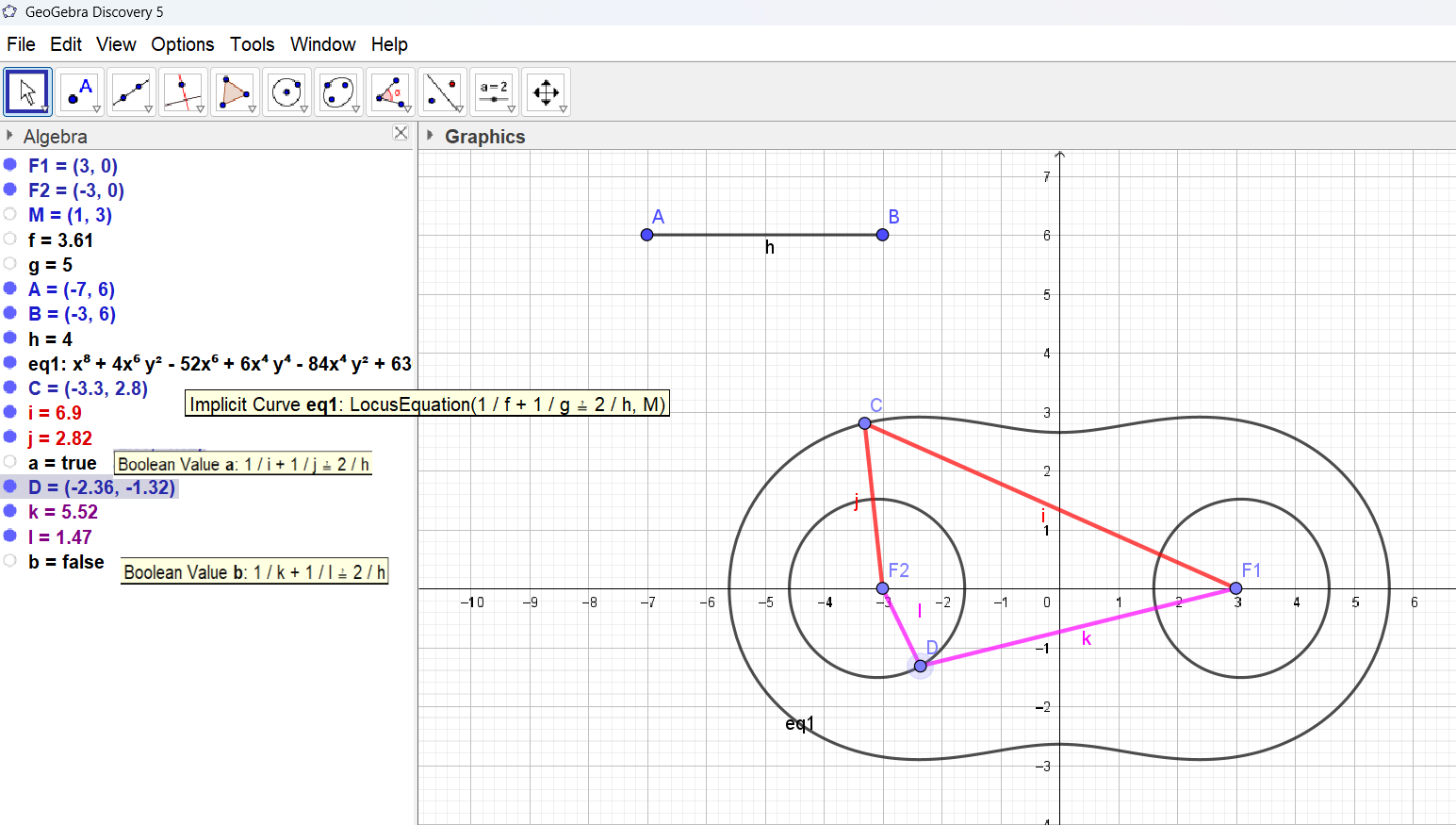}
\caption{Cayley ovals with superfluous components.}
\label{fig superfluous components}
\end{center}
\end{figure}

This is the kind of computation which is peformed behind the scene with GD's \textbf{LocusEquation}. The fact that the internal components are irrelevant can be assessed with GeoGebra-Discovery itself, as shown in Figure \ref{fig superfluous components}:  attach a point to the curve (GeoGebra has a command for this, and connect it by two segments to the foci $F_1$ and $F_2$. A screenshot of a corresponding session is displayed in Figure \ref{fig superfluous components}.
The bifocal definition from Definition \ref{def cayley oval} is verified for point $C$ and not for point $D$.

Now denote by $e=b/a$  (as an ellipse's eccentricity). Different shapes can be distinguished, according to the value of $e$.
\begin{enumerate}[(i)]
\item If $e<1$, the Cayley oval  has two components, namely two  non-intersecting disjoint loops, symmetric about the $y-$axis (Figure \ref{fig superfluous components}a).
\item If $e=1$, the curve looks like a lemniscate  (Figure \ref{fig superfluous components}b). Dynamic experiment shows how the external loops  of each pair in Figure \ref{fig superfluous components}a) tend to connect and form a lemniscate. 
\item If $1<e<\sqrt{3}$, the curve "looks like" a non-convex Cassini oval (Figure \ref{fig superfluous components}c). Of course, it is not a Cassini oval, we mean that at first glance we conjecture that the topology is similar, but the algebra says that the nature of the curve is different.
\item If $\sqrt{3}<e$, the curve is a (convex) oval, becoming more and more circular as $e$ increases (Figure \ref{fig superfluous components}d).
\end{enumerate}

 It is easy to check that equations obtained above with GD are special cases of the general polynomial equation (\ref{polyn eq cayley ovals}) computed wit Maple, for the specific values $b=2,3,4,6$ as indicated. All these polynomials are irreducible (proven with Maple, using the \textbf{evala(AFactor ($<$Polynomial$>$)} command, both on the general polynomial (\ref{polyn eq cayley ovals})) and on the special cases,   as in \cite{DPR 2023a}\footnote{It follows that the internal loops are not circles, but this is not our concern here.}

\section{Exploration of offsets of Cayley ovals}
\label{section offsets of Cayley ovals}
\begin{figure}[htb]
\begin{center}
 \includegraphics[width=4cm]{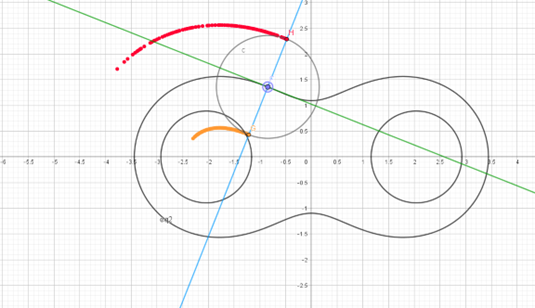}
\caption{The offset drawing technique.}
\label{fig drawing method}
\end{center}
\end{figure}

Since Cayley ovals may have four different shapes (some of them are not really ovals, as they have inflection points or more than one component), we explore each case separately. To plot an offset, we use the following method with the DGS (see Figure \ref{fig drawing method}). 
\begin{enumerate}[(i)]
\item Attach a point $F$ to the oval $\mathcal{C}$, so it will be possible to drag it only on $\mathcal{C}$.
\item Plot the tangent to $\mathcal{C}$ at $F$ using the \textbf{Tangent} command.
\item Plot the perpendicular to the tangent at $F$ using the \textbf{PerpendicularLine} command, i.e. plot the normal line to $\mathcal{C}$ at $F$.
\item Plot a circle centered at $F$ (by optional DGS button) with variable radius (use a slider or a variable segment which will be used as a slider when dragging one of its endpoints) and determine its points of intersection with the normal, using the \textbf{Intersect} command.
\item By activating \emph{Trace on} on the intersection points and dragging $F$ along $\mathcal{C}$, we obtain a preview of the offset.
\end{enumerate}

WLOG, we experiment with Cayley ovals, whose foci are on the $x-$axis, symmetric about the y-axis, and at a distance of 4 from each other. We started with small offset distance and increased it afterwards. In the following figures, we show some of the results.

\subsubsection{If $e<1$.}  
\begin{figure}[htb]
\begin{center}
 \includegraphics[width=4.5cm]{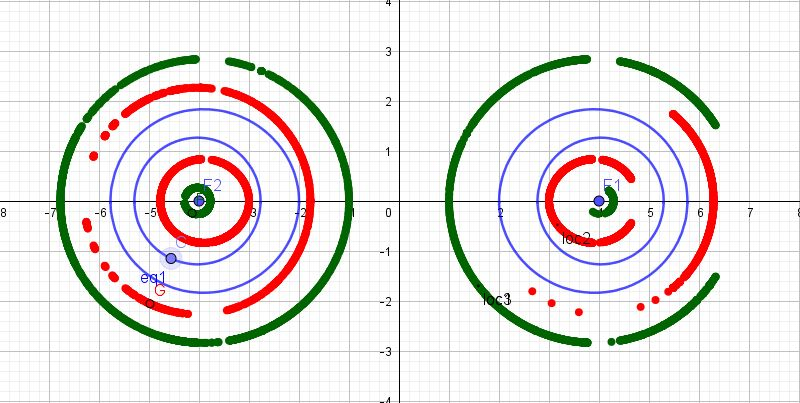}
\caption{Offsets of a Cayley oval with $e<1$}
\label{Experiments for offsets of Cayley ovals 1}
\end{center}
\end{figure}

Figure \ref{Experiments for offsets of Cayley ovals 1} shows a partial plot of the offset. Each loop provides two parallel loops.

\subsubsection{If $e=1$.} We obtain a unique shape of offsets: for small offset distances, the offset is the union of a loop (an oval) external to $\mathcal{C}$ and  2 internal loops  (Figure \ref{Experiments for offsets of Cayley ovals 2}) while for bigger offset distances, we obtain a rotated double-sided axe connected to 2 loops outside shaped-like curve that get growing as the offset increases (Figure \ref{Experiments for offsets of Cayley ovals 2}b,c)\footnote{An ongoing work with GeoGebra-Discovery developer Z. Kov\'acs (Linz, Austria) is aimed at improving the plots and eliminating the gaps in the plots}.

\begin{figure}[htb]
\begin{center}
 \includegraphics[width=3.5cm]{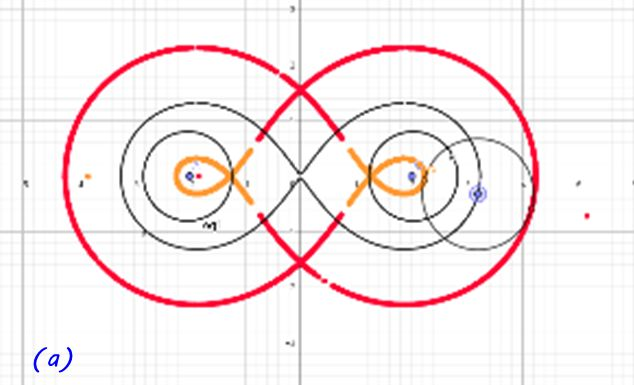}
 \qquad
  \includegraphics[width=3.5cm]{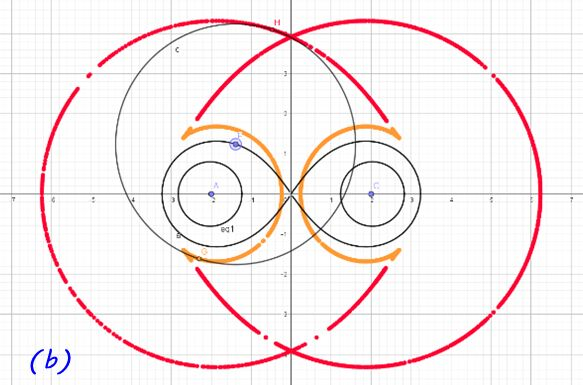}
 \qquad
  \includegraphics[width=3.5cm]{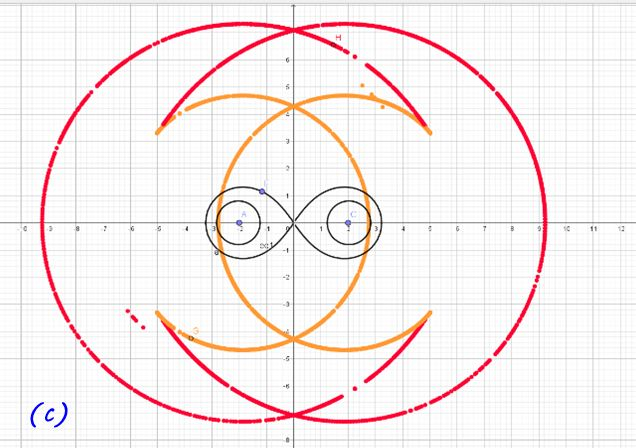}
 \caption{Offsets of a Cayley oval with $e=1$ at various distances}
\label{Experiments for offsets of Cayley ovals 2}
\end{center}
\end{figure}	

In each case, note the crunodes (i.e. the points of self-intersection). We study them in Section \ref{Singular points}.

\subsubsection{$1<e<\sqrt{3}$.} For small offset distances we get curves looking "really parallel" (Figure \ref{Experiments for offsets of Cayley ovals 2}a), but for some distance the inner and the outer curves are changed (not necessarily simultaneously). The inner curve became a three-loops shaped curve (Figure \ref{Experiments for offsets of Cayley ovals 2}b) and then changes twice to apple-like curves (Figure \ref{Experiments for offsets of Cayley ovals 2}a), and later to a sand-clock-like curves (Figure \ref{Experiments for offsets of Cayley ovals 2}d). The outer curve changes to a similar curve with an umbrella-like loops at the top and bottom of the offset (Figure \ref{Experiments for offsets of Cayley ovals 2}b,c,d).
\begin{figure}[htb]
\begin{center}
 \includegraphics[width=3.5cm]{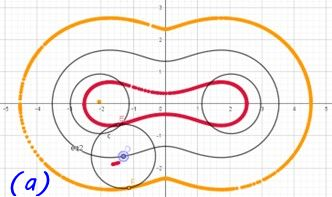}
 \qquad
  \includegraphics[width=3.5cm]{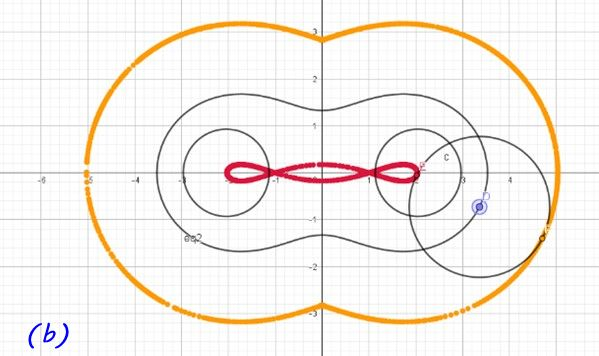}
 \qquad
  \includegraphics[width=3.5cm]{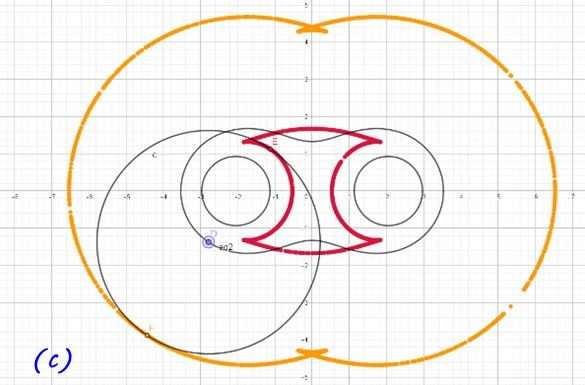}
 \qquad
  \includegraphics[width=3.5cm]{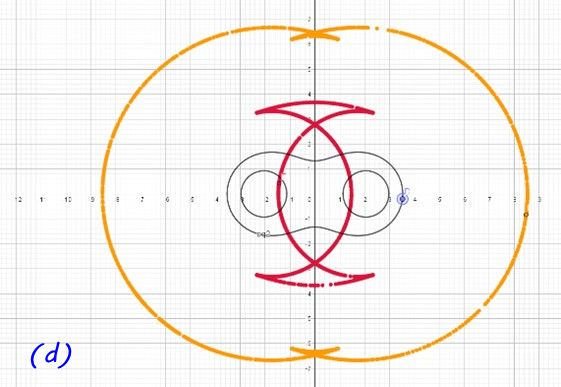}
 \caption{Offsets of a Cayley oval with $1<e<\sqrt{3}$ at various distances}
\label{Experiments for offsets of Cayley ovals 2}
\end{center}
\end{figure}	

\subsubsection{$e > \sqrt{3}$.}  
\begin{figure}[htb]
\begin{center}
 \includegraphics[width=3.5cm]{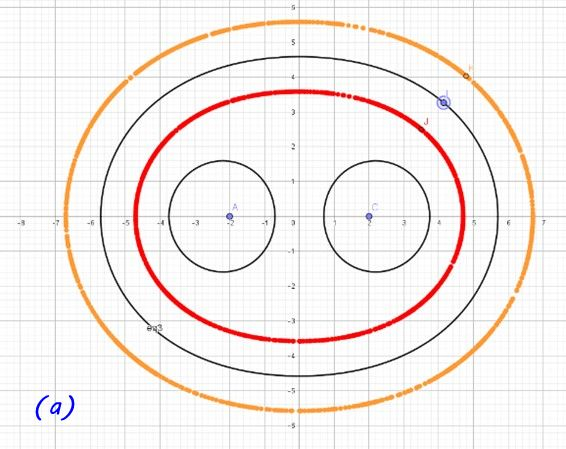}
 \qquad
  \includegraphics[width=3.5cm]{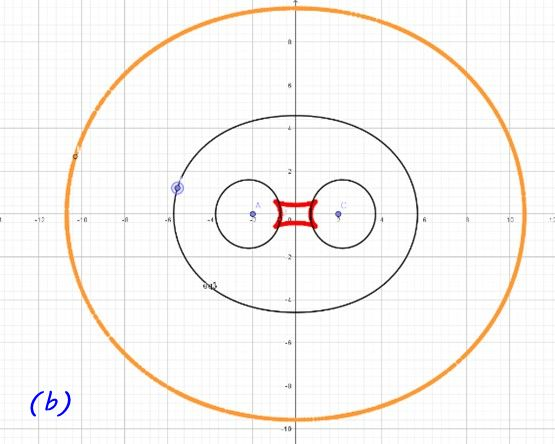}
 \qquad
  \includegraphics[width=3.5cm]{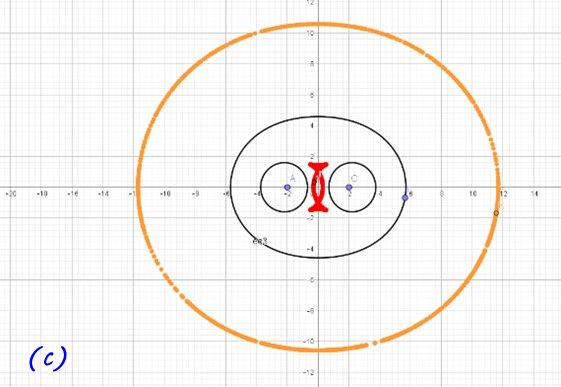}
 \qquad
  \includegraphics[width=3.5cm]{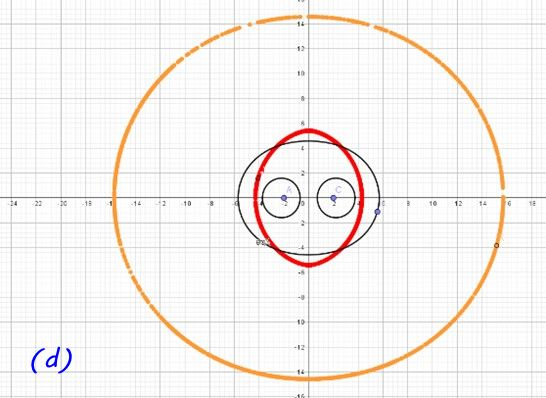}
 \caption{Offsets of a Cayley oval with $e>\sqrt{3}$ at various distances}
\label{Experiments for offsets of Cayley ovals 3}
\end{center}
\end{figure}

Here the ovals are dispatched into  2 types: the ellipse-like ovals and circle-like ovals.
For types, the external components are a parallel curves for every offset distance (Figures \ref{Experiments for offsets of Cayley ovals 3}). The internal components ares parallel curve for small offset distances (Figure \ref{Experiments for offsets of Cayley ovals 3}a), and changes gradually to rotated-astroid-like curve (Figure \ref{Experiments for offsets of Cayley ovals 3}b), then to a sand-clock-like curve (Figure \ref{Experiments for offsets of Cayley ovals 3}c) and afterwards to a kind of oval (Figure \ref{Experiments for offsets of Cayley ovals 3}d).

\section{Singular points}
\label{Singular points}

The singular points are dispatched into three categories (according to the taxonomy fixed by\cite{salmon}). \emph{Crunodes} are points of self-intersection, \emph{Acnodes} are isolated singular points and \emph{Cusps} are singular points where the curve has two different tangents. A special case exists for cusps, where the curve has a semi-tangent and "touches it form both sides". 
The above experimentations lead to studying crunodes and cusps.

\subsection{Mathematical analysis with software tools}
\label{subsection analysis with software}
The equations are too difficult to be solved by hand, therefore we work all along with software. In every case we display the code. As usual, we consider curves given by a parametrization $(x(t),y(t)), \; t \in \mathbb{R}$.
\begin{itemize}
\item Step 1: compute the first derivatives of $x(t)$ and $y(t)$
\item Step 2: determine  the outwards and inwards normals, denoted by: $nIn=(devY,-devx)$ and $nOut=(-devY,devX)$.  
\item Step 3: compute $(x(t),y(t))+d \; nIn$ and $(x(t),y(t))+d \; nOut$, where $d$ denotes the offset distance. we obtain the following parameterizations for the offsets:
\begin{itemize}
\item the internal offset curve parametrization:
 $\begin{cases} x=xI+ r \cdot x_{nIn} \\ y=yI+ r \cdot y_{nIn} \end{cases}$
\item the external offset curve parametrization:
 $\begin{cases} x=xO+ r \cdot x_{nOn} \\ y=yO+ r \cdot y_{nOn} \end{cases}$
\end{itemize}
\end{itemize}

Here is the Maple code; note that we run it for specific values of  $a$, $b$, and $r$ for the sake of visualization.

\footnotesize
\begin{verbatim}
with(plots);
x := b^2*cos(2*t)/(a*sin(2*t)^4);
y := sqrt(b^2/(4*cos(t)^4) - (x - a)^2);
dx := diff(x, t);
dy := diff(y, t);
v := Vector([dx, dy]);
len := sqrt(dx^2 + dy^2);
nI := Vector([dy/len, -dx/len]);
xI := r*nI[1] + x;
yI := r*nI[2] + y;
nO := Vector([-dy/len, dx/len]);
xO := r*nO[1] + x;
yO := r*nO[2] + y;
b := 1;
a := 1;
r := 1;
p1 := plot([xI, yI, t = 0 .. 20], color = "darkgreen");
p2 := plot([xO, -yO, t = 0 .. 15], color = "blue");
p3 := plot([x, y, t = 0 .. 20], color = "purple");
p4 := plot([x, -y, t = 0 .. 20], color = "purple");
display(p1, p2, p3, p4, scaling = constrained);
\end{verbatim}
\normalsize

\begin{figure}[htb]
\begin{center}
  \subfigure[d=0.1]{
 \includegraphics[width=3.5cm]{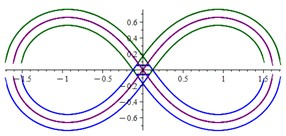}
}
 \qquad
  \subfigure[d=0.5]{
  \includegraphics[width=3.5cm]{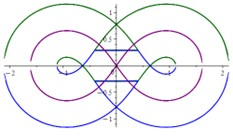}
}
 \qquad
  \subfigure[d=1]{
  \includegraphics[width=3.5cm]{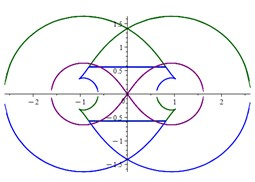}
}
 \qquad
  \subfigure[d=2]{
  \includegraphics[width=3.5cm]{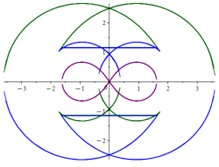}
}
 \qquad
  \subfigure[d=5]{
  \includegraphics[width=3.5cm]{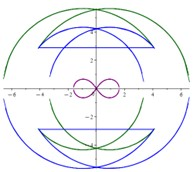} 
}
 \caption{Examples of offsets enhancing singular points}
\label{fig numerous offsets}
\end{center}
\end{figure}

The different offsets for $a=b=1$ and $d=0.1,0.5,1,2,5$ (in this order) are shown in Figure \ref{fig numerous offsets}.

Note that the plots show gaps instead of intersecting the x-axis. This is a frequent phenomenon with plots, in a neighborhood of a singular point. Moreover, superfluous horizontal segments. Adding standard options in the plot commands in order to avoid these segments led the computer to crash, after a long computation time. For us, this is still an open issue, that we will address in a further work. 
Moreover, we see, that the 2 colored plots of the  offsets are not clearly identified as to internal and external offsets, but we relate to them as positive and negative. That because the normal vector which has been defined flips whenever it crosses the $x-$axis (one more intriguing open issue). For that reason we plotted one offset regularly, and one using the symmetry about the $x-$axis. This appears in the change form $y$ to $-y$ in the code for $p3$ and $p4$.

The parametrization we used in the Maple session\footnote{It is given by \url{https://mathcurve.com/courbes2d.gb/cayleyovale/cayleyovale.shtml}} determines the external component of the Cayley oval. It does not determine the  "internal" components (the loops which appeared in Figure  \ref{Experiments for offsets of Cayley ovals 3}. Therefore, in Figure \ref{fig numerous offsets}, we see the offsets of the external loop of the Cayley oval only. Finding a parametrization for the internal loops is beyond the scope of the present paper. 

\subsection{First method: determination of the singular points according to the derivatives of $x(t)$ and $y(t)$.}
A singular point corresponds to a value $t_0$ of the parameter such that $(dx/dt)|_{t_0}$ and $(dy/dt)|_{t_0}$ either both vanish or at least one of them is not defined. we do not display the equations as they are heavy.
For the computations , we use Maple withe the command line called \textbf{solve}(devx=0,t). Actually, the command yields solutions over the field $\mathbb{C}$ of complex numbers, and the list has to be pruned in order to keep the real values only. This can be performed by using an additional option for the solve command. The coordinates of the singular points are then obtained by substitution of these values of the parameter in $xI yI$ or $xO yO$. 
Plotting the results gave us an interesting thing: the crunodes did not appear, although all of the cusps were found. 

As already mentioned, to find the cusps, we follow the next steps:
\begin{enumerate}[(i)]
\item Derive the x(t) offset curve parametrization 
\item Find the $t$ solutions whom fulfill the term $\frac{dx}{dt}=0$. For this purpose we will use the command solve in maple.
\item Filter the non-real values of $t$.
\end{enumerate}
Now we calculate the $x$ and $y$ for each value of $t$ we obtained in the previous section. We enhance the plots of the points to check accuracy.
Here is the Maple code used for this part of the work. Note we use here specific values for $a$,$b$, and $r$.

\footnotesize
\begin{verbatim}
with(plots):
x := b^2*cos(2*t)/(a*sin(2*t)^4);
y := sqrt(b^2/(4*cos(t)^4) - (x - a)^2);
dx := diff(x, t);
dy := diff(y, t);
v := Vector([dx, dy]);
len := sqrt(dx^2 + dy^2);
nI := Vector([dy/len, -dx/len]);
xI := r*nI[1] + x;
yI := -r*nI[2] - y;
nO := Vector([-dy/len, dx/len]);
xO := r*nO[1] + x;
yO := r*nO[2] + y;
b := 1;
a := 1;
r := 1;
dxi := diff(xI, t);
dxo := diff(xO, t);
sI := solve(dxi = 0, t);
sO := solve(dxo = 0, t);
t_valuesI := [];
for i to nops([sI]) do
    if Im(evalf(sI[i])) = 0 then t_valuesI := [op(t_valuesI), evalf(sI[i])]; 
    end if;
end do;
t_valuesO := [];
for i to nops([sO]) do
    if Im(evalf(sO[i])) = 0 then t_valuesO := [op(t_valuesO), evalf(sO[i])]; 
    end if;
end do;
pointsI := [seq([xI(t), yI(t)], t in t_valuesI)];
pointsO := [seq([xO(t), yO(t)], t in t_valuesO)];
curve_plotI := plot([xI, yI, t = 0 .. 7], color = darkgreen, thickness = 2);
points_plotI := pointplot(pointsI, symbol = solidcircle, symbolsize = 10, 
                color = red);
points_plotO := pointplot(pointsO, symbol = solidcircle, symbolsize = 10, 
                color = red);
curve_plotO := plot([xO, yO, t = 0 .. 7], color = blue, thickness = 2);
p1 := plot([x, y, t = 0 .. 20], color = "purple");
p2 := plot([x, -y, t = 0 .. 20], color = "purple");
display(curve_plotI, points_plotI, points_plotO, curve_plotO,p1,p2, 
        scaling = constrained);
\end{verbatim}
\normalsize

Now we run this code for different offset distances (its possible to change $a$ and $b$ as well); the output is in display in Figure \ref{fig numerous offsets 2} (for d=0.1,0.5,1,2,5 respectively):
\begin{figure}[htb]
\begin{center}
 \includegraphics[width=4cm]{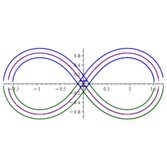}
 \qquad
  \includegraphics[width=4cm]{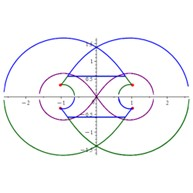}
 \qquad
  \includegraphics[width=4cm]{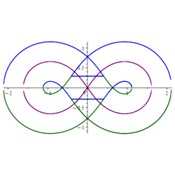}
 \qquad
  \includegraphics[width=4cm]{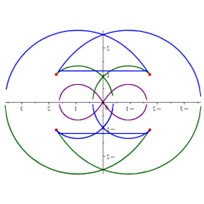}
 \qquad
  \includegraphics[width=4cm]{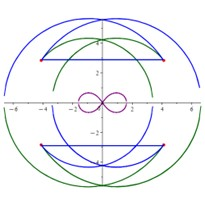} 
\caption{Examples of offsets at various distances enhancing singular points}
\label{fig numerous offsets 2}
\end{center}
\end{figure}	
The horizontal segments are irrelevant; they are a consequence of how the \textbf{plot} command works. Fine tuning is required to cancel them. All the pitfalls led to look for another approach, described in the next subsection. 

\subsection{Finding the singular points using the curvature}
Now, we try an another approach to get the cusp singular points. 	We use the method described in [25]. Let $\mathcal{C}$ be a curve determined by a parametric 
presentation $(x(t), y(t))$. If at a point $M(t) = (x(t), y(t))$, both functions $x$ and $y$ are 
differentiable at least twice, then the curvature of $\mathcal{C}$ at a this point is given by the formula: 
\begin{equation}
k(t)=\frac{r'\times r"}{v^3}  \cdot e_z=\frac{x'(t)y"(t)-x"(t)y'(t)}{(x'(t)^2+y'(t)^2)^{3/2}}
\end{equation}
 where $r$ is the distance of the offset and $v$ is the length of $r'$ . Define $\hat{k}=\frac{k}{\vert 1+kd \vert}$; solving the equation: k=$-\frac 1d$ yields the values of $t$ for which the curvature vanishes, i.e. the values of the parameter $t$ corresponding to the cusps. 
Here is the implementation in Maple code (note that we are initializing $r$ – for the sake of plotting):
\footnotesize
\begin{verbatim}
with(plots);
x := cos(t);
y := sin(t)^3;
d := 0.1;
dx := diff(x, t);
dy := diff(y, t);
denominator := sqrt(dx^2 + dy^2);
x_offset_pos := x + dy*d/denominator;
y_offset_pos := y - dx*d/denominator;
x_offset_neg := x - dy*d/denominator;
y_offset_neg := y + dx*d/denominator;
r := <x, y>;
dr := diff(r, t);
ddr := diff(dr, t);
kO := (-ddr[1]*dr[2] + ddr[2]*dr[1])/(dr[1]^2 + dr[2]^2)^(3/2);
kI := (ddr[1]*dr[2] - ddr[2]*dr[1])/(dr[1]^2 + dr[2]^2)^(3/2);
sO := solve(kO = -1/d, t);
sI := solve(kI = -1/d, t);
t_values := [];
for i to nops([sO]) do
    if Im(evalf(sO[i])) = 0 then t_values := [op(t_values), evalf(sO[i])]; exp(1)*nd; end if;
end do;
t_values1 := [];
for i to nops([sI]) do
    if Im(evalf(sI[i])) = 0 then t_values1 := [op(t_values1), evalf(sI[i])]; exp(1)*nd; end if;
end do;
outer_singular_points := [seq([evalf(subs(t = t_values[i], x_offset_pos)), 
    evalf(subs(t = t_values[i], y_offset_pos))], i = 1 .. nops(t_values))];
inner_singular_points := [seq([evalf(subs(t = t_values1[i], x_offset_neg)), 
    evalf(subs(t = t_values1[i], y_offset_neg))], i = 1 .. nops(t_values1))];
curve_plot := plot([x, y, t = 0 .. 2*Pi], color = purple);
outer_singular_points_plot := pointplot(outer_singular_points, 
    symbol = solidcircle, color = red, symbolsize = 8);
inner_singular_points_plot := pointplot(inner_singular_points, 
    symbol = solidcircle, color = red, symbolsize = 8);
offset_pos_plot := plot([x_offset_pos, y_offset_pos, t = -Pi .. Pi], 
    color = darkgreen, thickness = 2, discont = true);
offset_neg_plot := plot([x_offset_neg, y_offset_neg, t = -Pi .. Pi], 
    color = blue, thickness = 2, discont = true);
display(curve_plot, offset_pos_plot, offset_neg_plot, inner_singular_points_plot, 
    outer_singular_points_plot);
\end{verbatim}
\normalsize

\begin{figure}[htb]
\begin{center}
  \subfigure[d=0.75]{
  \includegraphics[width=4cm]{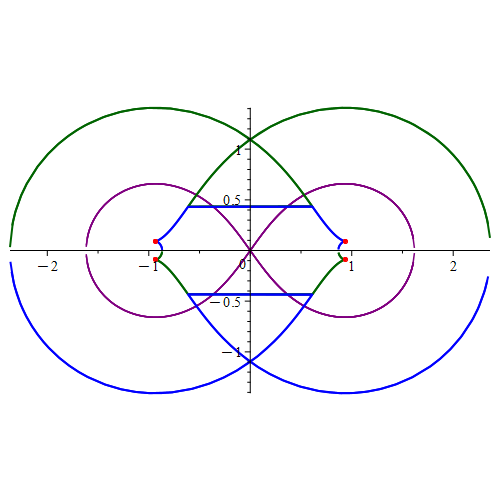}
}
 \qquad
 \subfigure[d=1]{
 \includegraphics[width=4cm]{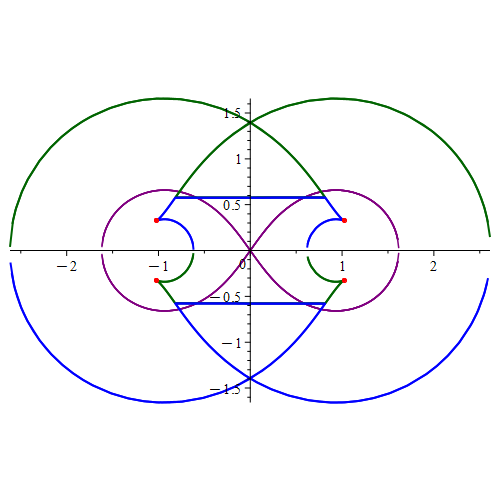}
}
 \qquad
  \subfigure[d=2]{
  \includegraphics[width=4cm]{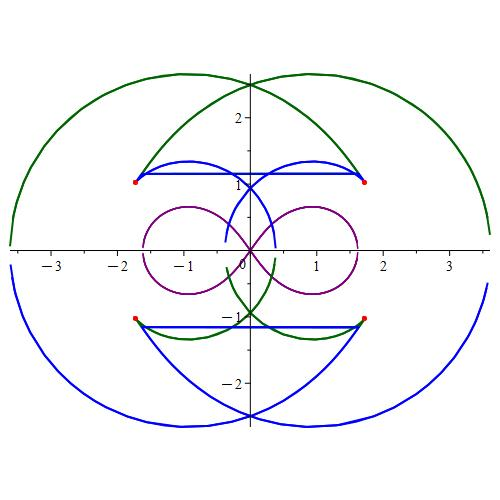}
}
 \qquad
\caption{Cusps which have been found according to the curvature}
\label{fig final plot}
\end{center}
\end{figure}

We run this code for the distances: d=0.75, 1, 2; and the output is displayed in Figure \ref{fig final plot} (in this order).

The first computations provide only half of the singular points. We solve this problem, using the symmetry about the $x-$axis, substituting  $-y$ instead of $-y$ as we did in subsection \ref{fig superfluous components}.

\subsection{Finding the crunodes}

Now, we will try to find the  points of self-intersection of Cayley ovals’ offsets. As mentioned – interestingly neither the method of $x(t)$ derivative nor the method of curvature give us the self-cross points. In this section we will try to find this type of singular points. 
We have to solve the following system of equations:
\begin{equation}
\label{solve eq for crunodes}
\begin{cases}
x(s)+\frac{\overset{\cdot}{y}(s) \cdot d}{\sqrt{\overset{\cdot}{x}(s)^2\overset{\cdot}{y}(s)^2}} =
x(t)+\frac{\overset{\cdot}{y}(t) \cdot d}{\sqrt{\overset{\cdot}{x}(t)^2\overset{\cdot}{y}(t)^2}}\\
y(s)+\frac{\overset{\cdot}{x}(s) \cdot d}{\sqrt{\overset{\cdot}{x}(s)^2\overset{\cdot}{y}(s)^2}} =
y(t)+\frac{\overset{\cdot}{x}(t) \cdot d}{\sqrt{\overset{\cdot}{x}(t)^2\overset{\cdot}{y}(t)^2}}
\end{cases}
\end{equation}

\begin{figure}[htb]
\begin{center}
  \subfigure[d=0.1]{
 \includegraphics[width=3.5cm]{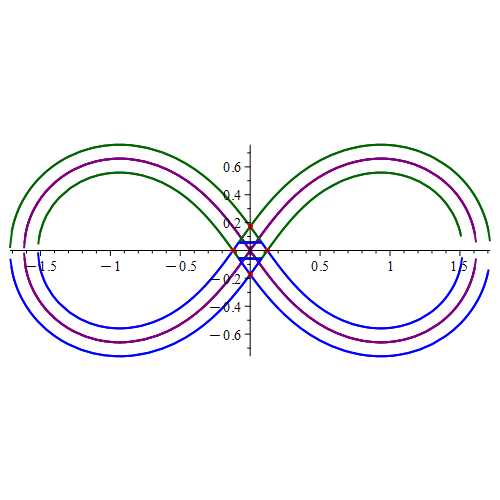}
}
 \qquad
  \subfigure[d=0.5]{
  \includegraphics[width=3.5cm]{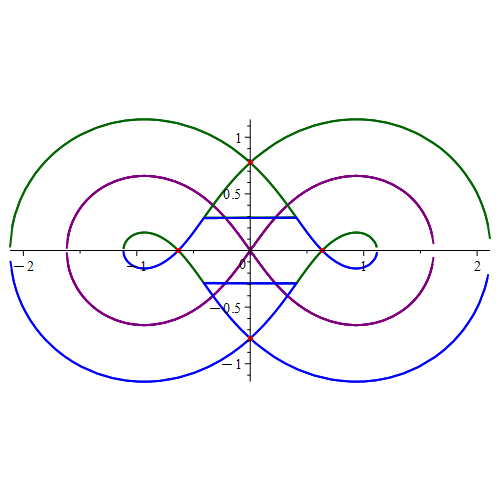}
}
 \qquad
  \subfigure[d=1]{
  \includegraphics[width=3.5cm]{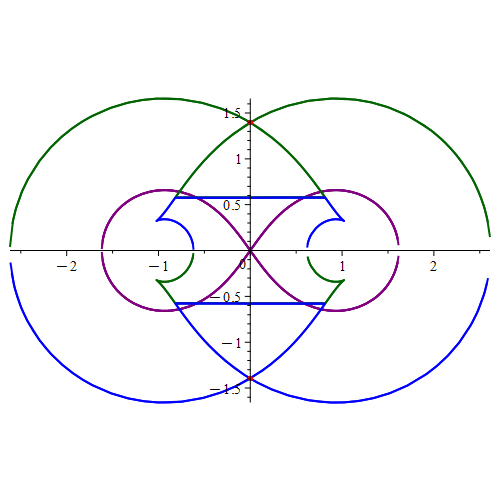}
}
 \qquad
  \subfigure[d=2]{
  \includegraphics[width=3.5cm]{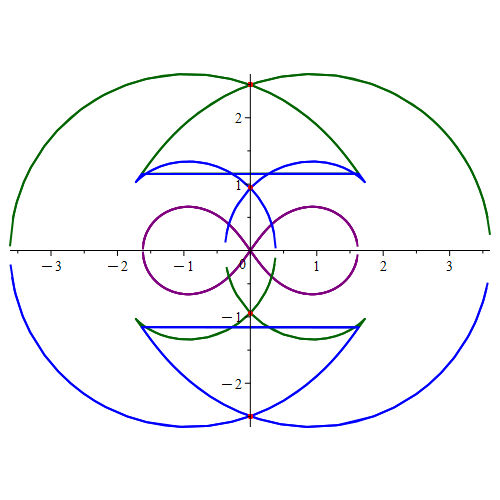}
}
 \qquad
  \subfigure[d=5]{
  \includegraphics[width=3.5cm]{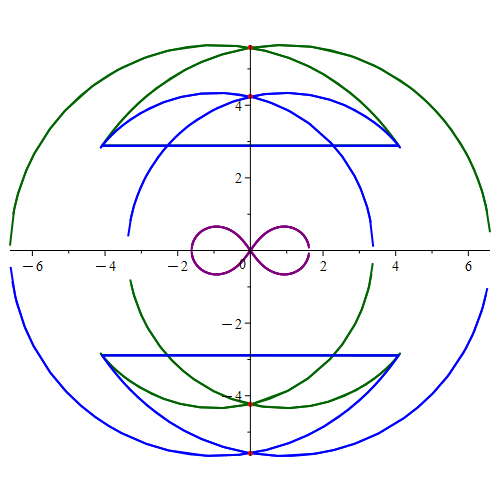} 
}
 \caption{Examples of offsets at various distances enhancing crunodes}
\label{fig numerous offsets 3}
\end{center}
\end{figure}

We implement it by the following code (Note we initializing $a$, $b$ and $r$ to specific values – for visualization reasons):
\footnotesize
\begin{verbatim}
with(plots);
solve_equations := proc(eqs, ranges)
	 local solutions, r, sol;
	 solutions := [];
	 for r in ranges do
		 sol := fsolve(eqs, r);
		 if sol <> NULL and evalf(subs(sol, s)) <> evalf(subs(sol, t)) then
		 	solutions := [op(solutions), sol];
	 	end if;
	 end do;
	 return solutions;
end proc;
extract_points := proc(solutions, x_func, y_func)
	 local points, sol, s_sol, t_sol;
	 points := [];
	 for sol in solutions do
		 s_sol := evalf(subs(sol[1], s));
		 t_sol := evalf(subs(sol[2], t));
		 if s_sol <> t_sol then
			 points := [op(points), [evalf(x_func(s_sol)), evalf(y_func(s_sol))]];
		 end if;
	 end do;
	 return points;
end proc;
a := 1;
b := 1;
d := 5;
x := u -> b^2*cos(2*u)/(a*sin(2*u)^4);
y := u -> sqrt(1/4*b^2/cos(u)^4 - (x(u) - a)^2);
neg_y := u -> -sqrt(1/4*b^2/cos(u)^4 - (x(u) - a)^2);
dx := D(x);
dy := D(y);
x_offset_pos := u -> x(u) + dy(u)*d/sqrt(dx(u)^2 + dy(u)^2);
y_offset_pos := u -> y(u) - dx(u)*d/sqrt(dx(u)^2 + dy(u)^2);
x_offset_neg := u -> x(u) - dy(u)*d/sqrt(dx(u)^2 + dy(u)^2);
y_offset_neg := u -> -y(u) - dx(u)*d/sqrt(dx(u)^2 + dy(u)^2);
eq1 := x_offset_neg(s) - x_offset_neg(t) = 0;
eq2 := y_offset_neg(s) - y_offset_neg(t) = 0;
eq3 := x_offset_pos(s) - x_offset_pos(t) = 0;
eq4 := y_offset_pos(s) - y_offset_pos(t) = 0;
eq5 := x_offset_pos(s) - x_offset_neg(t) = 0;
eq6 := y_offset_pos(s) - y_offset_neg(t) = 0;
ranges := [{s = 4*Pi .. 10*Pi, t = 4*Pi .. 10*Pi}, {s = -2*Pi .. 2*Pi, t = -2*Pi .. 2*Pi}];
solutions := solve_equations({eq1, eq2}, ranges);
solutions1 := solve_equations({eq3, eq4}, ranges);
solutions2 := solve_equations({eq5, eq6}, ranges);
points := [];
points := [op(points), op(extract_points(solutions, x_offset_neg, y_offset_neg))];
points := [op(points), op(extract_points(solutions1, x_offset_pos, y_offset_pos))];
points := [op(points), op(extract_points(solutions2, x_offset_pos, y_offset_pos))];
plot1 := plot([x(u), y(u), u = -Pi .. Pi], color = purple, thickness = 2, scaling = constrained);
plot2 := plot([x(u), -y(u), u = -Pi .. Pi], color = purple, thickness = 2, scaling = constrained);
plot3 := plot([x_offset_pos(u), y_offset_pos(u), u = -Pi .. Pi], color = darkgreen, thickness = 2,
	 scaling = constrained);
plot4 := plot([x_offset_neg(u), y_offset_neg(u), u = -Pi .. Pi], color = blue, thickness = 2,
	 scaling = constrained);
points_plot := pointplot(points, symbol = solidcircle, color = red, symbolsize = 7);
combined_plot := display(plot1, plot2, plot3, plot4, points_plot);
combined_plot;
\end{verbatim}
\normalsize
Now we run the code for $d=0.1,0.5,1,2,5$. The output is displayed in Figure \ref{fig numerous offsets 3}.

\begin{figure}[htb]
\begin{center}
 \includegraphics[width=3.5cm]{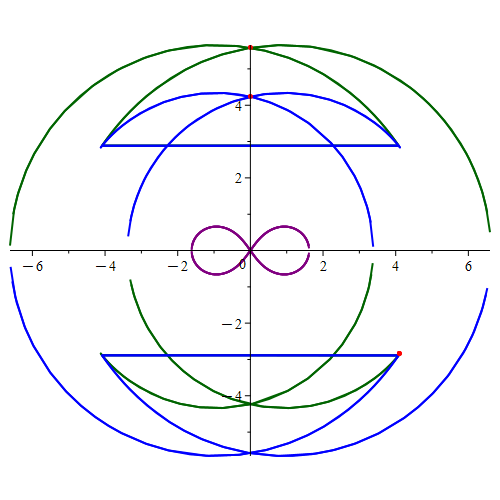}
 \qquad
  \includegraphics[width=3.5cm]{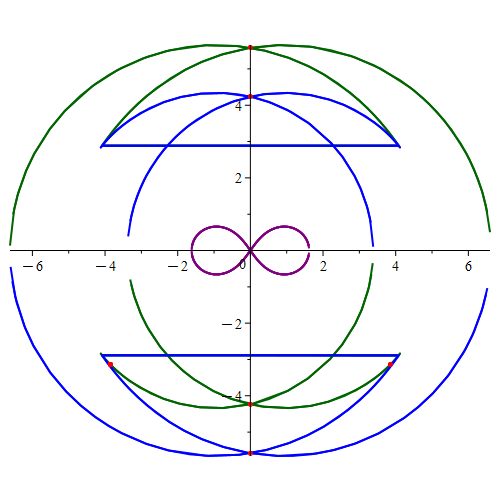}
\caption{Missing or unwanted points for wrong ranges when d=5}
\label{missing and unwanted crunodes}
\end{center}
\end{figure}

Note that we have solved every equation for several different intervals (called ranges in the code), because the algorithms provide generally  one solution in each given interval. Thus, we have to determine suitable specific range for every pair Cayley oval - offset distance for which we look for crunodes. For example, for $d=1$, 
$a=1$ and $b=1$, we computed the equations on the interval $[ 4\pi .. 10 \pi ]$, and not over the intervals shown in the code above. For that reason,  irrelevant  points may be obtained for certain intervals. Moreover, the opposite effect may occur and some points may be missed for a wrong choice of the interval (in Figure \ref{missing and unwanted crunodes}, some of the crunodes are enhanced, in red, and other ones do not appear).
Note also that the plots may be incomplete (e.g., in Figure \ref{missing and unwanted crunodes} the green and blue arcs should be connected, and they are not). This is a general issue; it received a nice solution with the \textbf{Plot2D} command ion GeoGebra-Discovery, but this command is not applicable here, because of a lack of an equation for the curve under study.

\section{Discussion}
\label{section discussion}
The study of offsets is non-intuitive. As already mentioned, the topology of an offset may be much more complicated than the topology of the progenitor. This appeared for offsets of an astroid  
\cite{offsets astroid} and of a kiss curve \cite{kiss curve}. In both cases, the progenitor has cusps, which induce specific properties on the offsets. In the present work, the progenitor has no singular point (with the only exception of a lemniscate), but anyway the offsets are quite complicated. There is no one-to-one correspondence between singularities of the progenitor and of the offset. Another wonderful example, with more simpler progenitor is provided by the 5-star in \cite{estrella}.

The usage of a Dynamic Geometry Software (DGS) enables to perform a rich exploration, leading to interesting discoveries and conjectures. These have to be proven, by  rigorous algebraic computations, made possible with a Computer Algebra System (CAS). Actually, the plot and animate features of the CAS complete the affordances of the DGS. Once again, the collaboration of a CAS and a DGS is fruitful.

When exploring these curves, Critical Thinking is crucial. The discussion in Section \ref{section Cayley ovals} addresses this issue. The 1st construction with a DGS of a Cayley oval, based on the bifocal definition of the oval, provided a curve with two irrelevant components. This has been pointed out using the DGS itself, its dynamics and its logical statements. Such a study is part of a wider filed in geometry, namely the study of geometric loci. Numerous studies have been published in the recent years, involving collaboration of CAS and DGS \cite{bp,hyperbolisms,pedals,jares-pech}; open problems still exist. This kind of activities ahem been successfully proposed by the 1st author to in-service teachers learning towards an advanced degree.

The offsets, which we presented here, have been studied per se. We wish to mention that offsets have numerous applications in engineering. Some of them have been briefly mentioned in \cite{DPK dialog}, for the design of industrial plants and entertainment parks. Others (offsets of Cassini ovals) appear in hydrology, for issues related to soil depollution \cite{wells}. 

Finally, we wish to mention that generative AI has not been used for this work. Asking several releases of software provided at most a general description of the topic, sometimes an outline of a possible way of working. Quite often, there were biases in the answer: even if the definition of an envelope or of an offset was correct, an example based n an ellipse gave already a wrong answer. For a comparison between an AI and GeoGebra in Geometry, see \cite{AI-ggb}. See also \cite{emprin-richard}. Anyway, with the state-of-the-art with CAS and DGS, the possible explorations are numerous and powerful. We wish to have a more efficient communication between kinds of software.  The sky is the frontier.

\paragraph{Acknowledgments:}
The first author acknowledges partial support from the CEMJ Chair at JCT. Many thanks to Zoltan Kov\'acs (Linz, Austria) for having accepted to help and mentor the 2nd author with DG.

\end{document}